\newcommand{\limfunc}[1]{\mathrm{#1}}
\newcommand{\func}[1]{\mathrm{#1}}

\documentclass[10pt,reqno]{amsart}
\makeatletter
\def\section{\@startsection{section}{1}%
	\z@{.7\linespacing\@plus\linespacing}{.5\linespacing}%
	{\bfseries%\normalfont\scshape
		\centering
}}
\def\@secnumfont{\bfseries}
\makeatother
\newtheorem{theorem}{Theorem}[section]

\newtheorem{proposition}[theorem]{Proposition}

\theoremstyle{definition}

\theoremstyle{remark}
\usepackage{graphicx}
\usepackage{pbox}

\usepackage{euscript}
\usepackage{amssymb}

\newtheorem{remark}[theorem]{Remark}
\numberwithin{equation}{section}
\usepackage{amsmath}
\usepackage{amsfonts}
\usepackage{amsthm}
\usepackage{array}
\usepackage{multirow}
\usepackage{graphicx}

\begin{document}
\title[Euler equations for Cosserat media]{Euler equations for Cosserat media}

\author[Valentin Lychagin]{Valentin Lychagin}
\thanks{The author was partially supported by the Russian Foundation for Basic Research (project 18-29-10013).}
\address{Valentin Lychagin: V. A. Trapeznikov Institute of Control Sciences of Russian Academy of Sciences, Moscow, 117997, Russia.}
\email{valentin.lychagin@uit.no}

\date{Date of Submission   August 15, 2020;  Date of Acceptance  October 10, 2020, Communicated by Yuri E. Gliklikh }

\begin{abstract}
We consider Cosserat media as $\mathbf{SO}(3)$-structures
over a domain $\mathbf{D}\subset\mathbb{R}^{3}$.
Motions of such media are given by
infinitesimal automorphisms of the $\mathbf{SO}(3)$-bundle.
We present Euler-type equations for such media and discuss
their structure.
\end{abstract}

\maketitle

%\noindent $\clubsuit$ Note to author: Use 2000 Mathematics Subject Classification.

\section{Introduction}

This paper is a realization of approach (\cite{DLT}) for the
case of media formed by `rigid microelements' (\cite{MGC}),
or Cosserat' media. 

We show that the geometry hidden behind of the such type of
media is a $\mathbf{SO}(3)$-structure over a spacial domain
$\mathbf{D}$. 

Namely, the dynamics in such media are given by vector
fields on the bundle $\pi\colon\Phi \rightarrow \mathbf{D}$
of $\mathbf{SO}(3)$-frames over $\mathbf{D}$. These fields
are solutions of the system of differential equations
(\cite{DLT}) that generalize the well known Navier-Stocks
equations.

To construct these equations we need two additional
geometric objects: (1) connection in the bundle $\pi$, we
call it media connection, and (2) left
$\mathbf{SO}(3)$-invariant metric on fibres of $\pi$. The
last determines the mechanics of `rigid microelements' and
the media connection allows us to compare `rigid
microelements' at different points of $\mathbf{D}$. The
third ingredient that we needed to construct the system of
differential equations is the thermodynamics of media. We
apply here the point of view (\cite{LY}), where
thermodynamics were considered as measurement of extensive
quantities.

Using these remarks and observations we present here  the
Euler type equations that govern the motion. To make the
presentation as compact as possible we restrict ourselves by
this type equations only although the difference between
Euler and Navier-Stocks type equations consists only in the
description of the stress tensor.

The final form of the Euler equations is given by
\eqref{EulerEq}, \eqref{MassConEq}, \eqref{ErgyConserv}.

\section{Geometry of $\mathbf{SO}(3)$}

Here we collect the main properties of the orthogonal
group $\mathbf{SO}(3)$.

Let $(\mathbf{T}, g)$ be a an Euclidian vector space, $\dim
\mathbf{T}=3$, where $g$ is a metric tensor. Let $A\in
\limfunc{End}\left( \mathbf{T}\right)$ be a linear operator
in $\mathbf{T}$ and let $A^{\prime}\in \limfunc{End}\left(
\mathbf{T}\right)$ be the $g$-adjoint operator, i.e
$g\left( AX,Y\right) =g\left(X,AY\right)$ for all $X,Y\in
\mathbf{T}$. Remind that 
\begin{equation*}
\mathbf{SO}(3)=
\left\{A\in\limfunc{End}
\left( \mathbf{T}\right)\, \vert\,
AA^{\prime}=\mathbf{1},\,\det A=1\right\} .
\end{equation*}
Geometrically, elements of the group are counterclockwise
rotations $R\left(\phi,\overline{n}\right)$
on angle $\phi$ about the axis through unit
vector $\overline{n}\in \mathbf{T}$. One has $R\left( \phi ,n\right)
=R\left( -\phi ,-n\right) $ and
$R\left( \pi ,n\right) =R\left( \pi,-n\right)$.
Therefore, as a smooth manifold, $\mathbf{SO}(3)$ is
diffeomorphic to the projective space $\mathbb{R}\mathbf{P}^{3}$.

The Lie algebra of the group, 
\begin{equation*}
\mathfrak{so}\left( 3\right) =\left\{ A\in \limfunc{End}\left( 
\mathbf{T}\right) \,\vert\, A+A^{\prime}=\mathbf{0}\right\},
\end{equation*}
consists of skew symmetric operators.

There is the \textit{hat isomorphism} of the Lie algebras
\begin{equation*}
\wedge\colon\left( \mathbf{T,\times }\right)
\rightarrow \mathfrak{so}\left(3\right),
\end{equation*}%
where $\left(\mathbf{T,\times}\right)$ is the Lie algebra
of vectors in $\mathbf{T}$ with respect
to the cross product $\times$.

In an orthonormal basis $\left( e_{1},e_{2},e_{3}\right)$
this isomorphism has the following form
\[
\wedge\colon
w=\left(w_{1},w_{2},w_{3}\right) \mapsto \widehat{w}=w_{1}
\widehat{e}_{1}+w_{2}\widehat{e_{2}}+w_{3}\widehat{e_{3}}=
\begin{bmatrix}
0 & -w_{3} & w_{2} \\ 
w_{3} & 0 & -w_{1} \\ 
-w_{2} & w_{1} & 0
\end{bmatrix}.
\]

\subsection{Exponent and logarithm}

In the case of Lie algebra $\mathfrak{so}\left( 3\right)$
the exponential map 
\begin{equation*}
\exp\colon\mathfrak{so}\left( 3\right) \rightarrow \mathbf{SO}(3),
\end{equation*}
has the following Rodrigues' form
(see, for example, \cite{Rod}, \cite{Mar}) 
\begin{equation*}
\exp \left( \phi \widehat{n}\right) =\mathbf{1+}\sin \left( \phi \right) 
\widehat{n}+(1-\cos \left( \phi \right) )\widehat{n}^{2},
\end{equation*}
and $\exp \left( \phi \widehat{n}\right) =R\left( \phi ,n\right)$.

This formula gives us the following description of the logarithm map
\begin{equation*}
\ln \left( R\right) =\frac{\arcsin \left( \phi \right) }{2\phi }\left(
R-R^{\prime }\right),
\end{equation*}
where $R\in \mathbf{SO}(3)$ and 
\begin{equation*}
\phi =\frac{3+\sqrt{-\limfunc{Tr}R^{2}}}{2}.
\end{equation*}

\subsection{Left invariant tensors on $\mathbf{SO}(3)$.}

The Baker--Campbell--Hausdorff formula in
$\mathfrak{so}\left(3\right)$ has very concrete form.

Namely, let $X,Y\in\mathfrak{so}\left(3\right)$, then
$\exp \left(X\right) \cdot \exp \left( Y\right) \in \mathbf{SO}(3)$ and
therefore has the form $\exp\left(Z(X,Y)\right)$,
for some element $Z(X,Y)\in \mathfrak{so}\left( 3\right)$.

Then (see, for example, \cite{Eng}),
\begin{equation*}
Z\left( X,Y\right) =\alpha\,X+\beta\,Y+\gamma\,\lbrack X,Y\rbrack,
\end{equation*}%
where%
\begin{equation*}
\alpha =\frac{a\arcsin \left( d\right) }{d\theta },\quad
\beta =\frac{b\arcsin\left( d\right) }{d\phi },\quad
\gamma =\frac{\arcsin \left( d\right) }{d\phi\theta },
\end{equation*}
$a$, $b$, $c$, $d$ are defined as follows
\begin{align*}
&a =\sin \left( \theta \right) \cos ^{2}\left( \frac{\phi }{2}\right)
-\omega \sin \left( \phi \right) \sin ^{2}\left( \frac{\theta }{2}\right),\\
&b=\sin \left( \phi \right) \cos ^{2}\left( \frac{\theta }{2}\right) -\omega
\sin \left( \theta \right) \sin ^{2}\left( \frac{\phi }{2}\right) \\
&c =\frac{1}{2}\sin \left( \theta \right) \sin \left( \phi \right) -2\omega
\sin ^{2}\left( \frac{\theta }{2}\right) \sin ^{2}\left( \frac{\phi }{2}
\right),\\
&d=\sqrt{a^{2}+b^{2}+2\omega ab+\left( 1-\omega ^{2}\right) c^{2}},
\end{align*}
and
\begin{equation*}
\theta =\sqrt{\frac{-\limfunc{Tr}\,X^{2}}{2}},\quad
\phi =\sqrt{\frac{-\limfunc{Tr}\,X^{2}}{2}},\quad
\omega =\theta^{-1}\phi^{-1}\sqrt{\frac{-\limfunc{Tr}\,XY}{2}}.
\end{equation*}%
Applying these formulae for the case $X=t\widehat{n},Y$,
where $n$ is a unit vector, we get
\begin{equation*}
\alpha =\frac{\phi }{2}\cot \left( \frac{\phi }{2}\right) +O\left( t\right)
,~\beta =1+t\left( \frac{\omega }{2}\cot \left( \frac{\phi }{2}\right) +%
\frac{\omega }{\phi }\right) +O\left( t^{2}\right) ,\ \gamma =\frac{1}{2}%
+O\left( t\right).
\end{equation*}

Denote by $E_{1}$, $E_{2}$, $E_{3}$ the left invariant vector fields on
$\mathbf{SO}(3)$ that correspond to the basis $e_{1}$, $e_{2}$, $e_{3}$
in $\mathfrak{so}\left( 3\right)$.

Then the above formulae give us the following expressions for
$E_{1}$, $E_{2}$, $E_{3}$ in the canonical coordinates
$\left(x_{1},x_{2},x_{3}\right)$ of the first kind:
\begin{align*}
E_{1}=&\frac{\phi }{2}\cot \left( \frac{\phi }{2}\right) \partial _{1}+
\frac{1}{2}\left( x_{2}\partial _{3}-x_{3}\partial _{2}\right) +\\
&x_{1}\left( 
\frac{1}{2}\cot \left( \frac{\phi }{2}\right) +\frac{1}{\phi }\right) \left(
x_{1}\partial _{1}+x_{2}\partial _{2}+x_{3}\partial _{3}\right),\\
E_{2}=&\frac{\phi }{2}\cot \left( \frac{\phi }{2}\right) \partial _{2}+
\frac{1}{2}\left( x_{3}\partial _{1}-x_{1}\partial _{3}\right) +\\
&x_{2}\left( 
\frac{1}{2}\cot \left( \frac{\phi }{2}\right) +\frac{1}{\phi }\right) \left(
x_{1}\partial _{1}+x_{2}\partial _{2}+x_{3}\partial _{3}\right),\\
E_{3}=&\frac{\phi }{2}\cot \left( \frac{\phi }{2}\right) \partial _{3}+
\frac{1}{2}\left( x_{1}\partial _{2}-x_{2}\partial _{1}\right) +\\
&x_{3}\left( 
\frac{1}{2}\cot \left( \frac{\phi }{2}\right) +\frac{1}{\phi }\right) \left(
x_{1}\partial _{1}+x_{2}\partial _{2}+x_{3}\partial _{3}\right),
\end{align*}%
where 
\begin{equation*}
\phi =\sqrt{x_{1}^{2}+x_{2}^{2}+x_{3}^{2}},\quad \text{and}\quad
\partial_{i}=\frac{\partial }{\partial x_{i}}.
\end{equation*}

Remark that basis vectors $e_{1},e_{2},e_{3}$ have the following commutation
relations:
$\lbrack e_{\sigma \left( 1\right) },e_{\sigma \left( 2\right) }\rbrack=
\limfunc{sign}\left( \sigma \right) e_{\sigma \left( 3\right) }$, for any
permutation $\sigma $ of three letters and accordingly vector fields
$E_{1}$, $E_{2}$, $E_{3}$ inherit the same commutation relations 
\begin{equation*}
\lbrack E_{\sigma \left( 1\right) },E_{\sigma \left( 2\right) }]=\limfunc{%
sign}\left( \sigma \right) E_{\sigma \left( 3\right) }.
\end{equation*}%
Let us denote by $\Omega _{1},\,\Omega_{2},\,\Omega_{3}\in
\Omega^{1}\left( 
\mathbf{SO}(3)\right)$ differential $1$-forms on Lie group
$\mathbf{SO}(3)$ such that $\Omega _{i}\left( E_{j}\right) =\delta _{ij},$
then 
\begin{equation*}
d\Omega _{\sigma \left( 3\right) }+\limfunc{sign}\left( \sigma \right) \
\Omega _{\sigma \left( 1\right) }\wedge \Omega _{\sigma \left( 2\right) }.
\end{equation*}%
Vector fields $E_{i}$ and differential $1$-forms $\Omega _{i}$ give us the
bases (over $\mathbb{R}$) in the space of left invariant vector fields and
correspondingly invariant differential $1$-forms on $\mathbf{SO}(3)$.

Moreover, any left invariant tensor on $\mathbf{SO}(3)$ is a linear
combination of tensor products $E_{i}$ and $\Omega_{j}$ with constant
coefficients.

Thus any left invariant metric $g$ on $\mathbf{SO}(3)$ is defined
by a positive self adjoint operator $\Lambda$ on
$\mathfrak{so}\left(3\right)$, so-called inertia tensor.

We will take basis $e_{1}$, $e_{2}$, $e_{3}$
to be eigenvectors of the operator $\Lambda$.

Thus we get:
\begin{equation}\label{Metr}
g_{\lambda }=\frac{1}{2}\left( \lambda _{1}\Omega _{1}^{2}+\lambda
_{2}\Omega _{2}^{2}+\lambda _{3}\Omega _{3}^{2}\right),
\end{equation}%
where $\lambda_{i}$, $\lambda _{i}>0$, are eigenvalues of
the operator and $\Omega _{i}^{2}$ are the symmetric squares
of the $1$-forms.

\subsection{Levi-Civita connections on $\mathbf{SO}(3)$}

Let $\nabla$ be the Levi-Civita connection associated with
left invariant metric $g_{\lambda}$. We denote by
$\nabla_{i}$ the covariant derivative along vector field
$E_{i}$.

Then we have 
\begin{equation*}
\nabla _{i}\left( E_{j}\right) =\sum_{k}\Gamma _{ij}^{k}E_{k},
\end{equation*}
where $\Gamma _{ij}^{k}$ are the Christoffel symbols.

This connection preserves the metric and therefore%
\begin{equation*}
g\left( \nabla _{i}\left( E_{j}\right) ,E_{k}\right) +g\left( E_{j},\nabla
_{i}\left( E_{k}\right) \right) =0,
\end{equation*}
or 
\begin{equation}\label{Christ1}
\lambda _{k}\Gamma _{ij}^{k}+\lambda _{j}\Gamma _{ik}^{j}=0,  
\end{equation}
for all $i,j,k=1,2,3.$

The condition for the connection to be torsion-free gives us
the following relations:
\begin{equation*}
\nabla _{\sigma \left( 1\right) }\left( E_{\sigma \left( 2\right) }\right)
-\nabla _{\sigma \left( 2\right) }\left( E_{\sigma \left( 1\right) }\right)
=E_{\sigma \left( 3\right) },
\end{equation*}%
for all permutations $\sigma$, or 
\begin{equation}\label{Christ2}
\Gamma _{\sigma \left( 1\right) ,\sigma \left( 2\right) }^{k}-\Gamma
_{\sigma \left( 2\right) ,\sigma \left( 1\right) }^{k}=\delta _{k,\sigma
\left( 3\right) }.
\end{equation}

The solution of these equations is the following%
\begin{eqnarray}\label{Gamma} 
\Gamma _{12}^{3} &=&\frac{\lambda -\lambda _{1}}{\lambda _{3}},\Gamma
_{23}^{1}=\frac{\lambda -\lambda _{2}}{\lambda _{1}},\Gamma _{31}^{2}=\frac{%
\lambda -\lambda _{3}}{\lambda _{2}},  \\
\Gamma _{21}^{3} &=&\Gamma _{12}^{3}-1,\Gamma _{32}^{1}=\Gamma
_{23}^{1}-1,\Gamma _{13}^{2}=\Gamma _{31}^{2}-1,  \notag
\end{eqnarray}%
where $\lambda =\lambda _{1}+\lambda _{2}+\lambda _{3},$
and all other Christoffel symbols are trivial.

Thus we have the only non-trivial relations:%
\begin{align*}
&\nabla_{1}\left( E_{2}\right)=\frac{\lambda -\lambda _{1}}{\lambda _{3}}%
[E_{1},E_{2}],\nabla _{2}\left( E_{3}\right) =\frac{\lambda -\lambda _{2}}{%
\lambda _{1}}[E_{2},E_{3}],\nabla _{3}\left( E_{1}\right) =\frac{\lambda
-\lambda _{3}}{\lambda _{2}}[E_{3},E_{1}],\\
&\nabla_{2}\left( E_{1}\right)=\frac{\lambda -\lambda _{2}}{\lambda _{3}}%
[E_{2},E_{1}],\nabla _{3}\left( E_{2}\right) =\frac{\lambda -\lambda _{3}}{%
\lambda _{1}}[E_{3},E_{2}],\nabla _{1}\left( E_{3}\right) =\frac{\lambda
-\lambda _{1}}{\lambda _{2}}[E_{1},E_{3}],
\end{align*}

\begin{theorem}
The Levi-Civita connection for left invariant the metric
$g_{\lambda}$ has the form:
\begin{align*}
\nabla_{\sigma \left( 1\right) }\left( E_{\sigma \left( 2\right) }\right) 
&=\limfunc{sign}\left( \sigma \right) \frac{\lambda -\lambda _{\sigma
\left( 1\right) }}{\lambda _{\sigma \left( 3\right) }}E_{\sigma \left(
3\right)},\\
\nabla_{i}\left( E_{i}\right)&=0,
\end{align*}
for all permutations $\sigma \in S_{3}$ and $i=1,2,3$.
\end{theorem}

\section{Cosserat media and $\mathbf{SO}\left(3\right)$-structures}

Let $\mathbf{D}$ be a domain in $\mathbb{R}^{3}$ considered
as Riemannian manifold equipped with the standard metric
$g_{0}$. Then by \textit{Cosserat medium} we mean a medium
composed with solids or having `rigid microstructure',
(\cite{MGC}). We assume that on the microlevel this media is
formed by rigid elements, which we represent as orthonormal
frames $f_{a}\colon\mathbb{R}^{3}\rightarrow
\mathbf{T}_{a}B$.

Thus configuration space of these media is the principal
$\mathbf{SO}\left(3\right)$-bundle
$\pi\colon\Phi\rightarrow\mathbf{D}$ of orthonormal frames
on $\mathbf{D}$, (see, for example, \cite{Ch}).
The projection $\pi$ assign
the centre mass $a$ of element $f_{a}$.

\subsection{Metrics and connections, associated with Cosserat medium}

The group $\mathbf{SO}\left( 3\right)$ acts in the natural
way on fibres of projection $\pi $  and we will continue to
use notations $E_{1}$, $E_{2}$, $E_{3}$ for the induced
vertical vector fields on $\Phi$. 

These fields form the basis in the module of
vertical vector fields on $\Phi$ and accordingly
differential $1$-forms
$\Omega_{1}, \Omega_{2}, \Omega_{3}$
define the dual basis in the space of vertical differential
forms.

We assume also that the media is characterized by
a $\mathbf{SO}\left(3\right)$-connection in the bundle $\pi$.
This connection, we call it 
\textit{media connection} and denote by $\nabla^{\mu}$, allows us to
compare microelements at different points of $\mathbf{D}$.

To define the connection, we consider a microelement as
orthonormal frame $b=\left( b_{1},b_{2},b_{3}\right)$,
formed by vector fields $b_{i}$ on $\mathbf{D}$.
Then the covariant derivatives $\nabla_{X}$ along vector field 
$X$ is defined by the connection form
$\omega\in\Omega^{1}\left( \mathbf{D,}\mathfrak{so}\left( 3\right) \right)$,
i.e. differential $1$-form on $\mathbf{D}$ with values
in the Lie algebra $\mathfrak{so}\left( 3\right)$,
and such that 
\begin{equation}
\nabla _{X}^{\mu }\left( b\right) =\omega \left( X\right) b.
\label{MediaConnection}
\end{equation}
By using the hat morphism, we represent the form as before
\begin{equation*}
\omega =%
\begin{Vmatrix}
0 & -\omega _{3} & \omega _{2} \\ 
\omega _{3} & 0 & -\omega _{1} \\ 
-\omega _{2} & \omega _{1} & 0
\end{Vmatrix},
\end{equation*}%
where $\omega_{i}$ are differential $1$-forms on $\mathbf{D}$.

Then formula \eqref{MediaConnection} shows us that
a microelement subject to rotation along vector
$\left(\omega_{1}(X),\,\omega_{2}(X),\,\omega_{3}(X) \right)$
on the angle
$\phi \left(X\right) =\sqrt{\omega_{1}\left( X\right) ^{2}+\omega _{2}\left( X\right)
^{2}+\omega _{3}\left( X\right) ^{2}}$, when we transport it on vector $X$.

Let $\left( x_{1},x_{2},x_{3}\right) $ be the standard Euclidian coordinates
on $\mathbf{D}$ and $\partial =\left( \partial _{1},\partial _{2},\partial
_{3}\right) $ be the corresponding frame. Here $\partial _{i}=\partial
/\partial x_{i}$.

In these coordinates we have 
\begin{equation*}
\omega _{i}=\sum_{j=1}^{3}\omega _{ij}d_{j},
\end{equation*}%
where $d_{j}=dx_{j}$.

The connection form $\omega$ will be the following
\begin{equation*}
\omega =
\begin{Vmatrix}
0 & -\omega _{31} & \omega _{21} \\ 
\omega _{31} & 0 & -\omega _{11} \\ 
-\omega _{21} & \omega _{11} & 0
\end{Vmatrix}
d_{1}+
\begin{Vmatrix}
0 & -\omega _{32} & \omega _{22} \\ 
\omega _{32} & 0 & -\omega _{12} \\ 
-\omega _{22} & \omega _{12} & 0
\end{Vmatrix}
d_{2}+
\begin{Vmatrix}
0 & -\omega _{33} & \omega _{23} \\ 
\omega _{33} & 0 & -\omega _{13} \\ 
-\omega _{23} & \omega _{13} & 0
\end{Vmatrix}
d_{3}.
\end{equation*}

This connection allows us to split tangent spaces
$T_{b}\Phi$ into direct sum 
\begin{equation*}
T_{b}\Phi =V_{b}+H_{b},
\end{equation*}
where $V_{b}$ is the vertical part with basis $E_{1,b}$,
$E_{2,b}$, $E_{3,b}$, and the horizontal space $H_{b}$ is
formed by vectors 
\begin{equation*}
X_{a}-\omega_{a}\left(X_{a}\right),
\end{equation*}
where $a=\pi\left(b\right)$.

Remark that geometrically spaces $H_{b}$ represent `constant
frames' due to \eqref{MediaConnection} and the basis in this
space is
\begin{eqnarray*}
&&\partial _{1}-\omega _{11}E_{1}-\omega _{21}E_{2}-\omega _{31}E_{3}, \\
&&\partial _{2}-\omega _{12}E_{1}-\omega _{22}E_{2}-\omega _{32}E_{3}, \\
&&\partial _{3}-\omega _{13}E_{1}-\omega _{23}E_{2}-\omega _{33}E_{3}.
\end{eqnarray*}%
The horizontal distribution
$H\colon b\in \Phi \rightarrow H_{b}\subset T_{b}\Phi$
could be also defined as the kernel of the following
differential $1$-forms on 
$\Phi$:
\begin{align*}
&\theta _{1} =\Omega _{1}+\omega _{11}d_{1}+\omega _{12}d_{2}+\omega
_{13}d_{3}, \\
&\theta _{2} =\Omega _{2}+\omega _{21}d_{1}+\omega _{22}d_{2}+\omega
_{23}d_{3}, \\
&\theta _{3} =\Omega _{2}+\omega _{31}d_{1}+\omega _{32}d_{2}+\omega
_{33}d_{3}.
\end{align*}

By metric associated with the left invariant metric
$g_{\lambda}$ and standard metric $g_{0}$ and media
connection form $\omega$ we mean direct sum of metric
$g_{\lambda}$ on the vertical space $V$ and the standard
metric $g_{0}$ on the horizontal space $H$.

We also call the media \textit{homogeneous} if the
connection form $\omega$ as well as the inertia tensor
$\Lambda$ are constants. Moreover, for the case of
homogeneous media the Euclidian coordinates
$\left(x_{1},x_{2},x_{3}\right)$ in the domain $\mathbf{D}$
are chosen in a such way that operators
$\widehat{\partial_{i}}\in \mathfrak{so}\left( 3\right)$
are eigenvectors of the inertia tensor $\Lambda$.

Summarizing we get the following
\begin{proposition}
The metric $g^{\mu}$ associated with the triple $\left( g_{\lambda
},g_{0},\omega \right)$ has the following form 
\begin{equation*}
g^{\mu }=\frac{1}{2}\sum_{i=1}^{3}\left( \lambda _{i}\Omega
_{1}^{2}+d_{i}^{2}\right).
\end{equation*}
\end{proposition}

\begin{remark}
Frame $\left( E_{1},E_{2},E_{3},\partial _{1}-\omega \left( \partial
_{1}\right) ,\partial _{2}-\omega \left( \partial _{2}\right) ,\partial
_{2}-\omega \left( \partial _{2}\right) \right)$ and coframe $\left( \Omega
_{1},\Omega _{2},\Omega _{3},d_{1},d_{2},d_{3}\right)$
are dual and $g^{\mu}$-orthogonal.
\end{remark}

\subsection{Levi-Civita connection on Cosserat media}

Assume that the media is homogeneous and let $\nabla$ be the
Levi-Civita connection on the configuration space $\Phi$
associated with metric $g^{\mu}$. Then for basic vector
fields $\partial_{i},E_{j}$, where $i,j=1,2,3$, we have the
following commutation relations 
\begin{equation*}
\lbrack\partial _{i},\partial_{j}\rbrack=
\lbrack\partial _{i},E_{j}\rbrack=0,\quad
\lbrack E_{\sigma(1)}, E_{\sigma(2)}\rbrack=
\limfunc{sign}\left( \sigma\right)E_{(3)}.
\end{equation*}

Moreover, it is easy to see that the equations for
Christoffel symbols do not depend on the connection form
$\omega$. On the other hand, it is clear that the Levi-Civita
connection $\nabla$, for the case $\omega =0$, is the direct
sum of the trivial connection on $\mathbf{D}$ and the left
invariant Levi-Chivita connection on $\mathbf{SO}\left(
3\right)$. Thus we get the following result.

\begin{theorem}
The Levi-Civita connection $\nabla^{c}$ on the configuration
space $\Phi$ associated with metric $g^{\mu}$ and
homogeneous media has the form, where the only non trivial
covariant derivatives are 
\begin{equation*}
\nabla _{E_{\sigma \left( 1\right) }}\left( E_{\sigma \left( 2\right)
}\right) =\limfunc{sign}\left( \sigma \right) \frac{\lambda -\lambda
_{\sigma \left( 1\right) }}{\lambda _{\sigma \left( 3\right) }}E_{\sigma
\left( 3\right) },
\end{equation*}
for all permutations $\sigma\in S_{3}$.
\end{theorem}

\subsection{Deformation tensor}

Dynamics we describe by $\pi$-projectable vector fields
on $\Phi$ (see \cite{DLT}, for more details):%
\begin{equation*}
U:=\sum_{i=1}^{3}\left( X_{i}\left( x\right) \partial _{i}+Y_{i}\left(
x,y\right) E_{i}\right).
\end{equation*}
Here $x=\left( x_{1},x_{2},x_{3}\right) $ are the Euclidian
coordinates on $\mathbf{D}$ and
$y=\left(y_{1},y_{2},y_{3}\right)$ are the canonical
coordinates on $\mathbf{SO}\left( 3\right)$.

Due to the above theorem and relations \eqref{Gamma} the covariant
differential $\mathcal{D}$ with respect to the Levi-Civita connection $%
\nabla ^{c}$ acts in the following way:%
\begin{eqnarray}
\mathcal{D}\left( \partial _{i}\right)  &=&0,\quad i=1,2,3,  \label{CovDif} \\
\mathcal{D}\left( E_{1}\right)  &=&\left( \alpha _{2}+1\right)\Omega
_{3}\otimes E_{2}+\left( \alpha _{3}-1\right) ~\Omega _{2}\otimes E_{3}, 
\notag \\
\mathcal{D}\left( E_{2}\right)  &=&\left( \alpha _{1}-1\right) \Omega
_{3}\otimes E_{1}+\alpha _{3}~\Omega _{1}\otimes E_{3},  \notag \\
\mathcal{D}\left( E_{3}\right)  &=&\alpha _{1}\Omega _{2}\otimes
E_{1}+\alpha _{2}\Omega _{1}\otimes E_{2},  \notag
\end{eqnarray}
where
\begin{align*}
\Gamma _{32}^{1}&=\alpha _{1}-1,\quad
& \Gamma _{23}^{1}&=\alpha _{1}=\frac{\lambda-\lambda _{2}}{\lambda _{1}}, \quad
& \Gamma _{31}^{2}&=\alpha _{2}+1,\\
\Gamma _{13}^{2}&=\alpha _{2}=\frac{\lambda_{1}}{\lambda _{2}}, \quad
& \Gamma _{21}^{3}&=\alpha _{3}-1,\quad
& \Gamma _{12}^{3}&=\alpha _{3}=\frac{\lambda-\lambda _{1}}{\lambda _{3}}.
\end{align*}%
Respectively, for the dual frame $\left( \Omega _{1},\Omega _{2},\Omega
_{3}\right) ,$ we have%
\begin{align}
\mathcal{D}\left( \Omega _{1}\right)  &=-\alpha _{1}\Omega _{2}\otimes
\Omega _{3}-\left( \alpha _{1}-1\right) \Omega _{3}\otimes \Omega _{2},
\notag \\
\mathcal{D}\left( \Omega _{2}\right)  &=-\alpha _{2}\Omega _{1}\otimes
\Omega _{3}-\left( \alpha _{2}+1\right) \Omega _{3}\otimes \Omega _{1}, 
\label{CovDifW} \\
\mathcal{D}\left( \Omega _{3}\right)  &=-\alpha _{3}\Omega _{1}\otimes
\Omega _{2}-\left( \alpha _{3}-1\right) \Omega _{2}\otimes \Omega _{1}. 
\notag
\end{align}

By the \textit{rate of deformation tensor}
$\Delta\left(U\right)$ we mean tensor 
\begin{multline*}
\Delta \left( U\right) =\mathcal{D}\left( U\right) =\sum_{i,j=1}^{3}\left(
\partial _{i}X_{i}\,d_{j}\otimes \partial _{i}+\partial
_{i}Y_{i}\,d_{j}\otimes E_{i}+E_{j}\left( Y_{i}\right) \Omega _{j}\otimes
E_{i}\right)+\sum_{i=1}^{3}Y_{i}\mathcal{D}\left( E_{i}\right) .
\end{multline*}

\section{Thermodynamics of Cosserat media}

The thermodynamics of the Cosserat media is based on
measurement (see, \cite{LY}, \cite{DLT}) of
\textit{extensive quantities:} \textit{inner energy}
$\mathit{E}$, \textit{volume} $V$, mass $m$, and deformation
$D$. The corresponding dual, or \textit{intensive
quantities} are the temperature $T$, preasure $p$, chemical
potential $\xi$ and the stress tensor $\sigma$.

The first law of thermodynamics requires that on
thermodynamical states the following differential $1$-form
\begin{equation}\label{thermo form1}
dE-(T\,dS-p\,dV+\limfunc{Tr}(\sigma^{\ast}dD)+\xi\,dm)  
\end{equation}
should be zero.

In other words, the thermodynamical state is a maximal
integral manifold $L$ of differential form \eqref{thermo form1}.

It is more convenient to use another but proportional
differential $1$-form
\begin{equation}\label{thermo form2}
dS-T^{-1}\left( dE+p\,dV-\limfunc{Tr}(\sigma ^{\ast }dD)-\xi\,dm\right).
\end{equation}

Here differential forms: $dE$, $T\,dS$,
$-p\,dV+\limfunc{Tr}(\sigma^{\ast }dD)$
and $\xi\,dm$ represent change of inner energy,heat,work and mass
respectively.

Extensivity of quantities $E,V,D,m,S$ means that their simultaneous
rescaling does not change the intensives as well as the thermodynamic state.

In other words, if we represent $\left(S, p, \sigma, \xi\right)$ as
functions of the extensive variables $\left(E,V,D,m\right)$,
then $S=S\left(E,V,D,m\right)$ is a homogeneous function of degree one.

Let $\left(s, \varepsilon, \rho, \Delta \right)$
be the densities of $\left(S, E, m, D\right)$.

Then substituting expression
$S=Vs(\varepsilon, \rho, \Delta)$ into \eqref{thermo form2}
we get that $1$-form 
\begin{equation*}
\psi=ds-T^{-1}\left( d\varepsilon -\limfunc{Tr}\,(\sigma ^{\ast }d\Delta
)-\xi\,d\rho \right) 
\end{equation*}%
should be zero and 
\begin{equation*}
\varepsilon -Ts=\limfunc{Tr}\,(\sigma ^{\ast }\Delta )-p+\xi \rho.
\end{equation*}%
The last condition shows that the density of Gibbs
free energy $\varepsilon-Ts$ equals sum of density
of deformation $\limfunc{Tr}\,(\sigma ^{\ast}\Delta )$
and mechanical $-p$ and chemical $\xi\,\rho$ works.

All these observations could be formulated as follows.

Let us introduce the thermodynamic phase space of Cosserat
medium (\cite{DLT}) as the contact space
$\widetilde{\Psi}=\mathbb{R}^{5}\times \limfunc{End}
\left( \mathbf{T}^{\ast }\right) \times \limfunc{End}
\left( \mathbf{T}\right)$, $\dim \widetilde{\Psi }=23$
with points $(s,T,\varepsilon,$ $\xi, \rho, \sigma,$
$\Delta) \in \Psi$,
where $\sigma \in \limfunc{End}\left( \mathbf{T}^{\ast }\right)$,
$\Delta \in \limfunc{End}\left( \mathbf{T}\right)$
and equipped with contact form $\psi$.
Then by thermodynamic
states we mean Legendrian manifolds
$\widetilde{L}\subset \widetilde{\Psi}$,
$\dim \widetilde{L}=11$, of the differential $\psi$, or respectively,
after eliminating entropy from the consideration,
their Lagrangian
projections $L\subset\Psi$ into
$\Psi=\mathbb{R}^{4}\times \limfunc{End}
\left( \mathbf{T}^{\ast }\right) \times
\limfunc{End}\left( \mathbf{T}\right)$, where the symplectic structure on $\Psi$ given by the
differential $2$-form 
\begin{equation*}
d\psi =\tau ^{-2}\left( dT\wedge d\varepsilon -\xi dT\wedge d\rho -dT\wedge 
\limfunc{Tr}\left( \sigma ^{\ast }d\Delta \right) \right) -T^{-1}\left( d\xi
\wedge d\rho -\limfunc{Tr}\left( d\sigma ^{\ast }d\Delta \right) \right).
\end{equation*}

As it was shown in (\cite{LY},\cite{DLT}) we should require, in addition,
that the differential symmetric form $\kappa $ shall define the Riemannian
structure on $L.$

Also, similar to (\cite{DLT}) we consider only such Legendrian manifolds $%
\widetilde{L}$ where $\left( T,\rho ,\Delta \right) $ are coordinates. In
this case we will write down form $T^{-1}\psi$ as  
\begin{equation*}
d(\varepsilon -Ts)-\left(s\,dT+\limfunc{Tr}\,\left( \sigma ^{\ast }d\Delta
\right) +\xi\,d\rho \right).
\end{equation*}%
Therefore, manifold $\widetilde{L}$ is Legendrian if and only if 
\begin{equation}\label{state eq}
s=h_{T},\quad\sigma =h_{\Delta },\quad\xi=h_{\rho},
\end{equation}%
where $h$ is the Gibbs free energy%
\begin{equation*}
h=\varepsilon -Ts.
\end{equation*}%
In the case of Newton-Cosserat media we have in addition
$\mathbf{SO}\left(3\right) \times \mathbf{SO}\left(3\right)$
symmetry and  Gibbs free energy $h$ is a function of
$\mathbf{SO}\left(3\right) \times \mathbf{SO}\left(3\right)$-invariants
(see \cite{DLT}, for more details).

Here we consider the Euler case (\cite{DLT}), when
\begin{equation*}
h\left(T,\rho,\Delta\right)=
p_{1}\left(T,\rho \right)\limfunc{Tr}\left( \Delta \right)+
p_{2}\limfunc{Tr}\left(\Delta\Pi_{V}\right),
\end{equation*}
where $\Pi_{V}$ is the projector on the vertical part of $T\Phi$.

In this case the stress tensor $\sigma$ equals 
\begin{equation*}
\sigma=
p_{1}\left(T,\rho\right)+
p_{2}\left(T,\rho\right)\Pi_{V},
\end{equation*}
or 
\begin{align}\label{EulerSigma}
\sigma=&\left( p_{1}\left( T,\rho \right) +p_{2}\left( T,\rho \right)
\right) ~\left( E_{1}\otimes \Omega _{1}+E_{2}\otimes \Omega
_{2}+E_{3}\otimes \Omega _{3}\right)+\\
&p_{2}\left( T,\rho \right) ~\left( \partial _{1}\otimes d_{1}+\partial
_{2}\otimes d_{2}+\partial _{3}\otimes d_{3}\right)
+p_{2}\sum_{i,j=1}^{3}\omega _{ij}E_{i}\otimes d_{j},\notag
\end{align}%
and energy density equals
\begin{equation*}
\varepsilon =\left( p_{1}-T p_{1,T}\right) \limfunc{Tr}\left( \Delta \right)
+\left(p_{2}-T p_{2,T}\right) \limfunc{Tr}\left(\Delta\Pi_{V}\right) 
\end{equation*}

\section{Euler equations for Cosserat media}

The general form of the Navier-Stocks equations and Euler equations as well,
for media with inner structures have the form (\cite{DLT}):

\begin{enumerate}
\item Moment conservation, or Navier-Stocks equation:
\begin{equation}\label{NSeq}
\rho \left( \frac{\partial U}{\partial t}+
\nabla_{U}^{c}\left( U\right)
\right) =\func{div}^{\flat }\sigma,
\end{equation}
where $\nabla^{c}$ is the Levi-Civita connection,
$\func{div}^{\flat}\sigma$ is a vector field dual to
the differential form $\func{div}\sigma\in \Omega^{1}
\left(\Phi\right)$, with respect to the canonical metric $g$.

\item Conservation of mass:
\begin{equation}\label{MassEq}
\frac{\partial \rho }{\partial t}+U\left( \rho \right) +\func{div}\left(
U\right) \rho =0.
\end{equation}

\item Conservation of energy:
\begin{equation*}
\frac{\partial \varepsilon }{\partial t}+
\varepsilon\,\func{div}\left( U\right)-
\func{div}\left( \zeta\,\func{grad}\left( T\right) \right) +
\limfunc{Tr}\left( \sigma ^{\prime }~\mathcal{D}U\right) =0,
\end{equation*}%
where $\zeta $ is the thermal conductivity.

\item State equations:%
\begin{equation*}
\sigma =h_{\Delta},\quad\varepsilon=h-T\,h_{T}.
\end{equation*}

\item In addition, we require that vector field
$U$ preserves the bundle
$\pi\colon\Phi\rightarrow\mathbf{D}$, or that $U$ is
a $\pi$-projectable vector field.
\end{enumerate}

In the case of Euler equations we have relation
\eqref{EulerSigma}, and using the property (see \cite{DLT}) 
\begin{equation*}
\func{div}\left( X\otimes \omega \right) =\func{div}\left( X\right) \omega
+\nabla _{X}\left( \omega \right),
\end{equation*}
where $X$ is a vector field and $\omega$ is
a differential $1$-form, we get, due to \eqref{CovDifW},
\begin{align*}
&\func{div}\left( E_{i}\otimes \Omega_{i}\right)=\nabla _{E_{i}}\left(
\Omega _{i}\right) =0,\\
&\func{div}\left( \partial _{i}\otimes d_{i}\right)=\nabla _{\partial
_{i}}\left( d_{i}\right) =0,
\end{align*}
and 
\begin{equation*}
\func{div}\,\sigma =\sum_{i=1}^{3}\left( E_{i}\left( p_{1}+p_{2}\right) \Omega
_{i}+(\partial _{i}-\sum_{j=1}^{3}\omega _{ji}E_{j})\left( p_{2}\right)
~d_{i}\right).
\end{equation*}%
Therefore,%
\begin{equation*}
\func{div}^{\flat}\,\sigma =\sum_{i=1}^{3}\left( \lambda _{i}^{-1}E_{i}\left(
p_{1}+p_{2}\right) E_{i}+(\partial _{i}-\sum_{j=1}^{3}\omega
_{ji}E_{j})\left( p_{2}\right) \partial_{i}\right).
\end{equation*}

Assume now that vector field $U$ has the form
\begin{equation*}
U=\sum_{i=1}^{3}\left( X_{i}\left( x\right) \partial _{i}+Y_{i}\left(
x,y\right) E_{i}\right).
\end{equation*}
Then,%
\begin{align*}
\nabla _{\partial _{1}}\left( U\right)  &=\sum_{i=1}^{3}\left( \partial
_{1}\left( X_{i}\right) \partial _{i}+\partial _{1}\left( Y_{i}\right)
E_{i}\right) ,\\
\nabla _{\partial _{2}}\left( U\right) &=\sum_{i=1}^{3}\left(
\partial _{2}\left( X_{i}\right) \partial _{i}+\partial _{2}\left(
Y_{i}\right) E_{i}\right) , \\
\nabla _{\partial _{3}}\left( U\right) &=\sum_{i=1}^{3}\left( \partial
_{3}\left( X_{i}\right) \partial _{i}+\partial _{3}\left( Y_{i}\right)
E_{i}\right) , \\
\nabla _{E_{1}}\left( U\right)  &=\sum_{i=1}^{3}\left( E_{1}\left(
Y_{i}\right) E_{i}+\lambda _{1i}Y_{i}[E_{1},E_{i}]\right)
=\sum_{i=1}^{3}E_{1}\left( Y_{i}\right) E_{i}+\lambda
_{12}Y_{2}E_{3}-\lambda _{13}Y_{3}E_{2}, \\
\nabla _{E_{2}}\left( U\right)  &=\sum_{i=1}^{3}\left( E_{2}\left(
Y_{i}\right) E_{i}+\lambda _{2i}Y_{i}[E_{2},E_{i}]\right)
=\sum_{i=1}^{3}E_{2}\left( Y_{i}\right) E_{i}-\lambda
_{21}Y_{1}E_{3}+\lambda _{23}Y_{3}E_{1}, \\
\nabla _{E_{3}}\left( U\right)  &=\sum_{i=1}^{3}\left( E_{3}\left(
Y_{i}\right) E_{i}+\lambda _{3i}Y_{i}[E_{3},E_{i}]\right)
=\sum_{i=1}^{3}E_{3}\left( Y_{i}\right) E_{i}+\lambda
_{31}Y_{1}E_{2}-\lambda _{32}Y_{2}E_{1},
\end{align*}
where
\begin{align*}
&\lambda_{12}=\frac{\lambda -\lambda_{1}}{\lambda _{3}},\quad
\lambda_{13}=\frac{\lambda -\lambda_{1}}{\lambda _{2}},\quad
\lambda_{21}=\frac{\lambda -\lambda_{2}}{\lambda _{3}},\\
&\lambda _{23}=\frac{\lambda -\lambda _{2}}{\lambda _{1}},\quad
\lambda _{31}=\frac{\lambda -\lambda _{3}}{\lambda _{2}},\quad
\lambda _{32}=\frac{\lambda -\lambda _{3}}{\lambda _{1}}.
\end{align*}

Therefore,
\begin{eqnarray*}
\nabla _{U}\left( U\right)  &=&\sum_{j,i=1}^{3}X_{j}\partial _{j}\left(
X_{i}\right) \partial _{i}+\sum_{j,i=1}^{3}X_{j}\partial _{j}\left(
Y_{i}\right) E_{i}+\sum_{j,i=1}^{3}Y_{j}E_{j}\left( Y_{i}\right) E_{i}+ \\
&&\left( \lambda _{23}-\lambda _{32}\right) Y_{2}Y_{3}E_{1}+\left( \lambda
_{31}-\lambda _{13}\right) Y_{1}Y_{3}E_{2}+\left( \lambda _{12}-\lambda
_{21}\right) Y_{1}Y_{2}E_{3}.
\end{eqnarray*}

Summarizing, we get the following system of Euler equations:
\begin{eqnarray}
&&\rho \left( \partial _{t}X_{i}+\sum_{j=1}^{3}X_{j}\partial _{j}\left(
X_{i}\right) \right) =(\partial _{i}-\sum_{j=1}^{3}\omega _{ji}E_{j})\left(
p_{2}\right) ,\ \ i=1,2,3,  \label{EulerEq} \\
&&\rho \left( \partial _{t}Y_{1}+\sum_{j=1}^{3}X_{j}\partial _{j}\left(
Y_{1}\right) +\sum_{j=1}^{3}Y_{j}E_{j}\left( Y_{1}\right) \right) +\left(
\lambda _{23}-\lambda _{32}\right) Y_{2}Y_{3}=\lambda _{1}^{-1}E_{1}\left(
p_{1}+p_{2}\right) ,  \notag \\
&&\rho \left( \partial _{t}Y_{2}+\sum_{j=1}^{3}X_{j}\partial _{j}\left(
Y_{2}\right) +\sum_{j=1}^{3}Y_{j}E_{j}\left( Y_{2}\right) \right) +\left(
\lambda _{31}-\lambda _{13}\right) Y_{1}Y_{3}=\lambda _{2}^{-1}E_{2}\left(
p_{1}+p_{2}\right) ,  \notag \\
&&\rho \left( \partial _{t}Y_{3}+\sum_{j=1}^{3}X_{j}\partial _{j}\left(
Y_{3}\right) +\sum_{j=1}^{3}Y_{j}E_{j}\left( Y_{3}\right) \right) +\left(
\lambda _{12}-\lambda _{21}\right) Y_{1}Y_{2}=\lambda _{3}^{-1}E_{3}\left(
p_{1}+p_{2}\right).  \notag
\end{eqnarray}%
Remark that the first three equations in \eqref{EulerEq}
are very close to the classical Euler equations and the
second three equations are the Euler type equations on
the Lie group $\mathbf{SO}\left(3\right)$.

The mass conservation equation takes the form
\begin{equation}\label{MassConEq}
\partial _{t}\rho +\sum_{i=1}^{3}\left( X_{i}\partial _{i}\rho
+Y_{i}E_{i}\rho \right) +\sum_{i=1}^{3}\left( \partial
_{i}X_{i}+E_{i}Y_{i}\right) \rho =0.  
\end{equation}

Finally, the energy conservation equation takes the form 
\begin{equation}
\frac{\partial \varepsilon }{\partial t}+\sum_{i=1}^{3}\left( \partial
_{i}X_{i}+E_{i}Y_{i}\right) ~\varepsilon +\sum_{i=1}^{3}E_{i}\left(
Y_{i}\right) =\func{div}\left(\zeta\,\func{grad}\left( T\right) \right) ,
\label{ErgyConserv}
\end{equation}
because 
\begin{equation*}
\limfunc{Tr}\left(\sigma^{\ast}\,\mathcal{D}U\right)
=\sum_{i=1}^{3}E_{i}\left( Y_{i}\right),
\end{equation*}%
and here%
\begin{equation*}
\varepsilon =\left( p_{2}-T\ p_{2,T}\right) \left( E_{1}\left( Y_{1}\right)
+E_{2}\left( Y_{2}\right) +E_{3}\left( Y_{3}\right) \right) +\left(
p_{1}-T~p_{1,T}\right) \func{div}\,U.
\end{equation*}


\begin{thebibliography}{99}
\bibitem{Rod} Rodrigues, Olinde, Des lois g\'{e}ometriques qui regissent les
d\'{e}placements d' un syst\'{e}me solide dans l' espace, et de la variation
des coordonn\'{e}es provenant de ces d\'{e}placement consid\'{e}r\'{e}es ind%
\'{e}pendant des causes qui peuvent les produire, J. Math. Pures Appl. 5
(1840), 380--440.

\bibitem{Eng} Kenth Eng\o , On the BCH-formula in so(3), BIT, Vol. 41, No.
3,(2001), pp. 629--632.

\bibitem{Mar} J. E. Marsden and T. S. Ratiu, Introduction to Mechanics and
Symmetry, Springer-Verlag, 1999

\bibitem{MGC} Holm Altenbach , G\'{e}rard A. Maugin, Vladimir Erofeev,
Mechanics of Generalized Continua, Springer-Verlag, pp. 371, 2011.

\bibitem{Ch} Chern, S. S.; Chen, W. H.; Lam, K. S. Lectures on differential
geometry. Series on University Mathematics, 1. World Scientific Publishing
Co., Inc., River Edge, NJ, 1999. x+356 pp

\bibitem{DLT} Duyunova Anna, Lychagin Valentin, Tychkov Sergey. Continuum
mechanics of media with inner structures. arXiv:2005.05840,2020, pp 14.

\bibitem{LY} Lychagin, Valentin, Contact geometry, measurement, and
thermodynamics. Nonlinear PDEs, their geometry, and applications, 3-52,
Tutor. Sch. Workshops Math. Sci., Springer, 2019.
\end{thebibliography}
\end{document}